\documentstyle[12pt]{article}
\newcommand{\argmin}{\mathop{\rm argmin}\limits}

\newcommand{\const}{\mathop{\rm const}\limits}

\newcommand{\mod}{\mathop{\rm mod}\limits}

\newcommand{\Ent}{\mathop{\rm Ent}\limits}

\begin{document}

\begin{center}

{\bf ADAPTIVE OPTIMAL SIGNAL (CARDIOGRAM) PROCESSING,\par
with boundary values and energy precisely measurement.} \\

\vspace{3mm}

   E.Ostrovsky, L.Sirota \ {\it Israel}. \\

\vspace{2mm}

{\it "O.L.E." Company; Rehovot, 76842, Israel; }\\
e \ - \ mails: galo@list.ru; \ eugostrovsky@list.ru; \\
\vspace{3mm}
{\it Department of Mathematics and Statistics, Bar-Ilan University,
59200, Ramat Gan.}\\
e \ - \ mail: \ sirota@zahav.net.il \\

\vspace{3mm}

                                   {\sc Abstract.}\\

\end{center}
\vspace{2mm}

  We construct an adaptive asymptotically  optimal in order in the {\it weight}
  Hilbert space norms signal denoising on the background noise and its energy measurement, with hight precision near the boundary of the signal.\par
  An offered method used the Fourier-Riesz expansion on the orthonormal polynomials,
for instance, Jacobi's polynomials, relative unbounded near the boundary weight function.\par
 An applications: technical and medical, in particular, cardiac diagnosis. \par

\vspace{2mm}

 \par

\vspace{3mm}

\section{Introduction. Statement of problem}
\vspace{3mm}

  A mathematical model: the cardiogram $ f = f(x),  x \in [-1,1]  $ is  Rieman integrable unknown function observed in the points $ \{ x(i) = x(i,n), i = 1,2,\ldots,n \} $ on the background additive white noise $ \epsilon(i), $  on the other terms:
{\sc regression problem (R): }

$$
\xi(i) = f(x(i)) + \sigma \cdot \epsilon(i), \eqno(0)
$$
where $ \{ \epsilon(i) \} $ is an additive noise;
$  \epsilon(i) $ are independent centered: $ {\bf E} \epsilon(i)=0 $ normed random variables (r.v.), $ {\bf Var} [\epsilon(i)] = 1, $ with unknown distribution; in particular, r.v. $\epsilon(i) $ not necessary to be Gaussian distributed; so that

 $$
 \sigma^2 = {\bf Var} [\xi(i)] = \const \in (0,\infty),
 $$
$$
x(i) = x(i,n) = (2i-n-1)/(n-1),  \  (n \ge 15),
$$
$ f = f(x) $ be a Rieman integrable unknown function (on the other terms, signal,
for instance, cardiogram).\par

{\bf Our aim is to offer the so-called adaptive optimal in order as $ n \to \infty$
in the modular Hilbert space sense estimation of a function (signal) $ f =  f(x) $ and its energy of arbitrary order with a greet precision especially near the boundaries $ x = \pm 1. $ }\par
\vspace{3mm}
 In the ordinary Hilbert spaces this problem is considered in many publications; see, for instance, \cite{Bobrov1}, \cite{Cai1},  \cite{Donoho1}, \cite{Efromovich1}, \cite{Klemela1}, \cite{Lepsky1}, \cite{Polyak1}, \cite{Tchentsov1} and reference therein.\par
\vspace{3mm}

\section{Notations and definitions}

\vspace{3mm}

{\bf 1.Modular Hilbert spaces.} Let $ X = [-1,1], x \in X, \Phi = \{ \phi_k(x) \},
\ k = 1,2,\ldots $ be some fixed complete bounded:
$$
\sup_k \sup_{x \in [-1,1]} |\phi_k(x)| < \infty,
$$
uniform Lipshitz continuous:

 $$
 |\phi_k(x) - \phi_k(y)| \le C \cdot k \cdot |x - y|
 $$
 orthonormal polynomial system relative some non-trivial non-negative integrable
 {\it weight signal} function  $ \gamma = \gamma(x): $

$$
(\phi_k,\phi_l)_{\gamma} := \int_{-1}^{1} \phi_k(x) \phi_l(x)  \gamma(x) dx =
\delta(k,l),
$$
$ \delta(k,l) $ is the Kroneker's symbol;
 $$
 \deg{\phi_k} = k-1, \ k = 1,2,\ldots.
 $$

{\it We impose henceforth the following conditions on the weight function} $ \gamma(x): \
 \gamma(x) \ge 0, $

$$
K(\gamma) \stackrel{def}{=} (2\pi)^{-1} \int_{-1}^1 \frac{\gamma(x)}{\sqrt{1-x^2}} \ dx=(2\pi)^{-1} \int_{-\pi}^{\pi} \gamma(\cos \phi) d \phi < \infty; \eqno(1.a)
$$

$$
\int_{-1}^1 \frac{\log(\gamma(x))}{\sqrt{1-x^2}} \ dx > - \infty; \eqno(1.b)
$$

$$
\exists \kappa = \const < 1/2, \ \forall x \in (-1,1) \ \Rightarrow
|\gamma^/(x)| \le C \left[(x+1)^{-\kappa-1} + (x-1)^{-\kappa-1}  \right].\eqno(1.c)
$$

 We redefine the weight function $ \gamma = \gamma(x) $ at the boundary points
 $ x = \pm 1 $ as follows:
 $$
 \gamma(-1) = \gamma(+1) = 0.
 $$
 The conditions (1.a), (1.b) and (1.c) are satisfied, e.g. for the weight function

 $$
 \gamma_{\alpha,\beta} (x) = (1-x)^{\alpha} (1+x)^{\beta}, \ \alpha,\beta > -1/2, \
 x \in (-1,1);
 \gamma_{\alpha,\beta} (\pm 1) =0,
 $$
 for the classical (normed!) Jacobi's  polynomial, which we will denote

 $$
 \phi^{(\alpha,\beta)}_k(x): \ \deg(\phi^{(\alpha,\beta)}_k(x)) = k-1,
$$
 $$
 \int_{-1}^1 \phi^{(\alpha,\beta)}_k(x) \phi^{(\alpha,\beta)}_l(x)
 \gamma_{\alpha,\beta} (x) dx = \delta(k,l).
 $$
 {\it As long as we need to estimate with greet precision the behavior of unknown
 function $ f(x) $ near the points $ x = \pm 1, $ we will suppose further that}

$$
\alpha,\beta \in (-1/2, 0]. \eqno(2)
$$

 Note that
$$
 K(\gamma_{\alpha,\beta} ) = 2^{\alpha + \beta} \ \pi^{-1} \ B(\alpha+ 1/2,\beta+ 1/2),
$$
obviously when $ \alpha,\beta > -1/2; \ B(\alpha,\beta)  $ denotes ordinary Beta- function. \par
 Many used for us information about orthogonal polynomials see in \cite{Geronimus1},
\cite{Simon1}, \cite{Szego1}. For instance, there are described the recurrence relations
of a view

$$
\phi_{k+2}(x) = (A(k;\Phi)x + B(k; \Phi)) \phi_{k+1}(x) +
C(k; \Phi)\phi_k(x), k=0,1,2,\ldots,
$$
 with initial conditions
$$
 \phi_0(x) = \left[ \int_{-1}^1 \gamma(z) \ dz \right]^{-1/2} =: C_0 = C_0(\gamma),
$$

$$
 \phi_1(x) = \lambda(x - C_1),
$$
where

$$
C_1 = C_1(\gamma) = \frac{\int_{-1}^1 x \gamma(x) dx}{\int_{-1}^1 \gamma(x) dx},
$$

$$
\lambda^{-2} = \lambda^{-2}(\gamma) = \int_{-1}^1 (x-C_1)^2 \ \gamma(x) \ dx.
$$

 For instance,
$$
\int_{-1}^1 \gamma_{\alpha,\beta}(x) dx = 2^{\alpha+\beta+1}  \ B(\alpha+1,\beta+1),
$$

$$
C_1(\gamma_{\alpha,\beta})= \frac{\beta-\alpha}{\alpha+\beta+2},
$$

$$
\lambda^{-2}(\gamma_{\alpha,\beta}) = 2^{\alpha+\beta+1} B(\alpha+1,\beta+1) \times
$$

$$
\frac{3\alpha^2\beta+3\alpha \beta^2+3\alpha^2 + 3\beta^2 + 16\alpha \beta + 14\alpha
+14\beta + 12}{(\alpha+\beta+3)(\alpha+\beta+2)^2}.
$$

These recurrent relations with correspondent initial conditions
 may be used for the programming of offered algorithms.\par

{\bf Remark 1.} By our opinion, it is convenient for the practical using to choose
the interval (0,1) instead (-1,1) and the weight function

$$
\gamma(x) = \frac{x^{\alpha-1} (1-x)^{\beta-1}}{B(\alpha,\beta)}, \ x \in (0,1);
$$
$ \gamma(0)= \gamma(1)= 0; \ \alpha,\beta > 0; $  i.e. $ \gamma(x) $ is the density of
the well-known Beta distribution. \par
Then $ \phi_0(x) = 1, $

$$
\phi_1(x) = \frac{(\alpha+\beta)x -\alpha}{\sqrt{\alpha \beta/(\alpha+\beta+1)}}.
$$

 The useful values of the parameters $ \alpha,\beta  $ are following:

 $$
 \alpha,\beta \in (1/2,1].
 $$

\vspace{4mm}
  For arbitrary integrable function $ f:[-1,1] \to R $ we put
 $$
f(x) = \sum_{k=1}^{\infty} c(k)\phi_k(x);
$$
here
$$
 c(k) = \int_{-1}^1 f(x) \ \phi_k(x) \ \gamma(x) \ dx
$$
 be the Fourier-Riesz coefficients of a function $ f(\cdot) $
 over the system $ \Phi. $ \par

  Let $ w = w(k), k=1,2,\ldots $  be some positive function (sequence), for example,
  $ w(k) = k^{m}, \ m \ge 0. $ The value of $ m $ is said to be the order of energy;
  not necessary to be integer. \par
  We  will called the function $ w=w(k) $ the {\it energy weight, } in contradiction to the notion {\it signal weight} $ \gamma = \gamma(x). $\par

 We introduce the following Hilbert modular space $ B = B(\Phi,w), $ which will called {\it modular Hilbert weight space,} consisting on
all the (measurable) functions $ \{f \} $ with finite quadratic functional

 $$
W(w)[f] = W(\Phi,w)[f] \stackrel{def}{=}
 \sum_{k=1}^{\infty} w(k)|c(k)|^2.\eqno(3)
$$
 The Hilbert's norm in this space may be defined as usually:
 $$
 ||f||W(w) = \sqrt{W(w)[f]}.
 $$
  For instance, $ W(w)[f] $ may be integral on the square of {\it derivative} of an 
  arbitrary order of the signal $ f(\cdot): $ 
  
  $$
  W(w)[f] = \int_{-1}^1 [f^{(r)}(t)]^2 \ dt
  $$
  or its weight modification:
  $$
  W(w)[f] = \int_{-1}^1 \zeta(t) [f^{(r)}(t)]^2 \ dt
  $$
    
 In the case when $ w(k) = w(k) = w_{\theta}(k) \asymp k^{\theta}, \
 \theta = \const \ge 0$ we will denote $ W(\Phi, w_{\theta} ) = W_{\theta}(\Phi)
  = W_{\theta}. $ \par
 Evidently, if $ w(k) = w_0(k) = 1, $  then

 $$
 \sum_{k=1}^{\infty} c^2_k = ||f||^2(\Phi, w_0) = (f,f)_{\gamma}.
 $$

 For instance, the Sobolev's spaces $ W_2^m = W_2^m[-1,1], m = 0,1,2,\ldots $
consisting on all the periodical $ (\mod 2)  \ m-$ times differentiable functions $ f $ with finite norm

$$
||f||^2W_2^m = \left[(||f||L_2[-1,1])^2 + (||f^{(m)}||L_2[-1,1])^2 \right]^{1/2}
$$
is weighted Hilbert modular space relative the trigonometric system $ T. $ \par
 The functional

$$
W(f) = W(f; \Phi, w)=  ||f||^2W(\Phi,w)
$$
will be called {\it Energy } or more exactly $ w \ - $ {\it Energy } of the signal $ f. $ We denote any estimation of $ W(f) $ as $ W(n, f). $\par

\vspace{3mm}

{\bf 2. Loss functions.} We define the following loss function
$ L = L(g(n,\cdot), f(\cdot)) $ for an arbitrary estimation of the function $ f $
(more exactly, the sequence of estimations)  $ \ g(n,x) = g(n,x; \{\xi(i)\}) $ as
 $$
 L_f(g(n,\cdot), f(\cdot))= L^{(\Phi)}_{\gamma}(g(n,\cdot), f(\cdot)) \stackrel{def}{=}
 {\bf E} ||g(n,\cdot)-f(\cdot)||^2_{\gamma}. \eqno(4.a)
 $$

We denote as for the arbitrary energy estimation $ W(n,f) $ as $ M(W(n,f), W(f) )$  the following loss function (functional):

$$
Z(W(n,f), W(f) ) \stackrel{def}{=} {\bf E} (W(n,f) - W(f))^2.\eqno(4.b)
$$

\vspace{3mm}

{\bf 3. Denotations.} We denote:

$$
S(N;m) = \sum_{k=1}^N w^m(k), \ S(N) = S(N,1) = \sum_{k=1}^N w(k);
$$

$$
\rho_w(N;f)= \rho_w(N) = \sum_{k=N+1}^{\infty} w(k) \ |c(k)|^2, \
\rho(N) = \sum_{k=N+1}^{\infty} |c(k)|^2, \eqno(5)
$$
then
$$
 f \in  W(\Phi,w) \Leftrightarrow \lim_{N \to \infty}\rho(N;f,1) = 0.
$$

\vspace{3mm}

\section{ Adaptive signal estimation: denoising}

\vspace{3mm}

{\bf 1. Offered estimations.}  Let us denote the {\it truncated} Fourier-Riesz's sum for the function $ f(x) $  as follows:

$$
f_{(N)}(x) = \sum_{k=1}^N c(k) \phi_k(x),
$$
and correspondingly

$$
f_{(N,n)}(x) = \sum_{k=1}^N c_n(k) \phi_k(x),
$$
where

 $$
 c_n(k) = n^{-1}\sum_{i=1}^n f(x(i,n)) \cdot \phi_k(x(i,n)) \cdot
 \gamma(x(i,n)).
 $$
 The {\it non-random} variables  $ c_n(k) $ are usually consistent  approximations of a Fourier (or Fourier-Riesz) coefficients $ c(k); $  they are Rieman integral sums for
 the function $ f(x)\phi_k(x) \gamma(x). $ \par

 Let us consider a so-called projection, or Tchentzov's {\it estimation} $f(N,n,x)$
 of the function  $  f(\cdot) $ of a view

$$
f(N,n,x) = \sum_{k=1}^N c(k,n) \phi_k(x). \eqno(6a)
$$
where {\it  random variables }
$ c(k,n) $ are usually consistent  estimations of a Fourier coefficients $ c(k) $  based on the date  $ \xi(i), i = 1,2,\ldots,n, $ namely, for the considered regression problem  {\sc R}

$$
c(k,n) = n^{-1}\sum_{i=1}^n \xi(i) \cdot \phi_k(x(i,n)) \cdot \gamma(x(i,n)).\eqno(6b)
$$

We denote also
$$
A(N,n) = \rho(N) + \frac{N}{n}.
$$
Notice that

$$
A(N,n) \asymp  {\bf E}||f(N,n,\cdot) - f(\cdot)||^2_{\gamma}, \ N = 1,2,\ldots, [n/3];
$$
$$
A^*(n): = \min_{N \in [1,n/3]}A(N,n); \ N_0: = \argmin_{N \in [1,n/3]} A(N,n).\eqno(7)
$$

{\it We assume in addition that}

$$
\lim_{N \to \infty} \frac{N_0}{\sqrt{n}} = 0. \eqno(7a)
$$

  This condition guarantee us that

  $$
\lim_{n \to \infty} \frac{ |  f_{N_0,n} - f_{N_0}(x)|^2_{\gamma}}{A^*(n)} = 0.
\eqno(7b)
  $$

 For example, if $ \gamma(x) = 1 $ and  $ \rho(N) \asymp N^{-2\Delta},
 N \to \infty,  $ then the last condition implies  the well-known  inequality $ \Delta > 1/2, $  see, e.g. \cite{Juditsky1}, \cite{Lepsky1}. \par
 Analogous conclusion is true when $ \rho(N) \asymp N^{-2\Delta} L(N), $  where
 $ L(N) $ is slowly varying as $ N \to \infty $ function.\par
  The function $ N \to N^{-2\Delta}L(N) $ is called regular varying as $ N \to \infty $ 
function. Used for us facts about slowly and regular varying function see in the 
classical monograph of E.Seneta \cite{Seneta1}.  \par
  The condition (7a) holds also when $ \rho(N) \asymp N^{-2\Delta}L(N), \ \Delta > 1/2, $
 and if the signal weight function $ \gamma(x) $ satisfied the conditions (1a),(1b)
 and (1c).\par

\vspace{3mm}

  Further, we introduce the following important (random) statistics (functionals):
$$
\tau(N,n) = \sum_{k=N+1}^{2N} |c(k,n)|^2; \ M = M(n) =
\argmin_{N \in [1,n/3]} \tau(N,n), \eqno(8)
$$

$$
\tau^*(n) = \min_{N \in [1,n/3]} \tau(N,n) = \tau(M(n),n) \eqno(9)
$$
and we introduce {\it our estimation} $ \hat{f}(n,x) $ of a function $ f(x) $ as follows:
$$
\hat{f}(n,x) = f(M(n),n,x) = \sum_{k=1}^{M(n)} c(k,n) \phi_k(x). \eqno(10)
$$
 Note that $ M(n) $ is a random sequence which depended only on the source sample
  $ \{ \xi(i) \}, i = 1,2,\ldots,n $ (adaptivity). \par

\vspace{3mm}

{\bf 2.} {\it We suppose also}

$$
0 < \underline{\lim}_{N \to \infty} \rho(2N)/\rho(N) \le \overline{\lim}_{N \to \infty} \rho(2N)/\rho(N) < 1. \eqno(11a)
$$
Another case, i.e. when, e.g.,

 $$
 \rho(2N)/\rho(N) > 0, \  \lim_{N \to \infty} \rho(2N)/\rho(N) = 0 \eqno(12)
  $$
may be investigated analogously. \par

{\bf 3.} Obviously, in the considered here regression problem {\sc R } we must impose
some conditions on the errors of measurements $ \epsilon(i). $ {\it We suppose in
addition that the centered and normed variables} $ \{ \epsilon(i) \} $  {\it satisfy the following condition: there exist a finite positive constants} $ q,Q $ {\it such that for all positive values} $ u $

  $$
  \max \left[ {\bf P}(\epsilon(i) > u),{\bf P}(\epsilon(i) < -u) \right] \le
  \exp \left(-  (u/Q)^q \right). \eqno(13a)
  $$
 It may be assumed in many practical cases that the r.v.$ \{\epsilon(i)\} $ are Gaussian 
 distributed; in this  case the inequality (13a) holds with parameters $ q=2 $ and 
 $ Q=1. $\par

\vspace{3mm}
\section{ Main result}
\vspace{3mm}

If we choose $ N = N_0, $ then we conclude that the optimal speed of convergence for the {\it non-adaptive estimation} is asymptotically $ A^*(n): $ as
$ n \to \infty $
$$
{\bf E}||f(n,N_0,\cdot)-f(\cdot)||^2 \asymp A^*(n). \eqno(14)
$$
 At the same result is true for our {\it adaptive estimation.} \par
 Recall that we suppose the conditions (1a), (1b), (1c), (7a), (11a) and (13a) to be
 satisfied.\par
 {\bf Theorem 1.} Let $ K $ be arbitrary regular function space with a norm
 $ ||\cdot||K $ compactly embedded in the space $ B = B(\Phi,w).$ We denote the minimax value of the loss function as $ A^*(n)= A^*(n,K,Z), Z = \const \in (0,\infty): $
$$
A^*(n,K,Z)\stackrel{def}{=} \inf_{g(n,\cdot)}\sup_{f \in K, ||f||K \le Z} L(g(n,\cdot), f(\cdot)), \eqno(15)
$$
where  $ "\inf" $ is calculated over all, {\it including even  non adaptive} estimations of $ f(\cdot).$ \par
 We assert that under conditions for adaptive estimation
our estimation $ \hat{f}(n,\cdot) $ is asymptotically as $ n \to \infty $ optimal up to finite multiplicative constant in the minimax sense on the arbitrary regular function classes $ K $ compactly embedded in the modular space $ W = W(\Phi,w)$:

$$
\sup_{f \in K, ||f||K \le Z} L(\hat{f}(n,\cdot), f(\cdot)) \le C(\gamma,\Phi,K,Z) \cdot A^*(n,K,Z), \ C(\gamma,\Phi,K,Z) < \infty. \eqno(16)
$$

{\bf 5. Variance estimation.}
Note that the consistent $ 1/\sqrt{n} $ estimation $ \sigma^2(n) $ of the value
$ \sigma^2 $ may be constructed by formulae \cite{Lepsky1}, \cite{Polyak1}

$$
\sigma^2(n) = (n - M(n) - 1)^{-1}
\sum_{i=1}^n \left( \hat{f}(n,x(i,n))- \xi(i)  \right)^2.\eqno(17)
$$

\vspace{3mm}

{\bf 6. Confidence region.} For the rough building of the confidence domain, also
 adaptive, in the $ B(\Phi,\gamma) $ norm we proved the following result. \par
{\bf Theorem 2.} We assert under at the same conditions as in theorem 1 that there
exist positive finite {\it non-random} constructive and estimative  constants
$ Y_j = Y_j(\gamma, \Phi,K,Z), \ j=1,2 $  and for the problem (R) $ Y_j =
Y_j(\gamma,\Phi,K,Z;q,Q) $ such that

$$
\overline{\lim}_{n \to \infty} \tau^*(n)/A^*(n,K,Z)= Y_1 \eqno(18)
$$
and correspondingly
$$
\overline{\lim}_{n \to \infty} || \hat{f}(n,\cdot) - f(\cdot))||^2_{\gamma}/\tau^*(n)= Y_2. \eqno(19)
$$

Therefore, we conclude: with probability tending to one as $ n \to \infty $

$$
|| \hat{f}(n,\cdot) - f(\cdot)||^2_{\gamma} \le 1.05 \cdot Y_2 \cdot \tau^*(n). \eqno(20)
$$

 More exact computation show us that

 $$
 {\bf E} ||\hat{f}(n,\cdot) - f(\cdot)||^2_{\gamma} \sim {\bf E}\tau^*(n),
 \eqno(21)
 $$

 $$
 {\bf Var} ||\hat{f}(n,\cdot) - f(\cdot)||^2_{\gamma} \sim
 4 \sigma^2 K^2(\gamma) \rho(M(n))/n + 3 \sigma^4(n) K^4(\gamma) M(n)/n^2
 \stackrel{def}{=} \delta^2(n), \eqno(22)
  $$
where

$$
\tau^*(n) \sim \rho(M(n)) + \sigma^2(n) K^2(\gamma) M(n)/n.\eqno(23)
$$
 Therefore, with probability tending to 0.95

 $$
 ||\hat{f}(n,\cdot) - f(\cdot)||^2_{\gamma} \le \tau^*(n) + 2.54 \cdot \delta(n).
 \eqno(24)
 $$

\vspace{3mm}

\section{Energy estimation}

\vspace{3mm}
 In order to detect a distortion in the signal, for instance, infarction, heart attack, apoplexy in the cardiodiagnostics; leakage in the pipelines, in particular, in the nuclear reactors  we offer to estimate its energy and comprise with it
 standard value up to confidence  domain. \par
 {\bf 1. Ordinary energy.}\par
 Let us define  for simplicity the so-called  ordinary energy  functional

 $$
 G = G(f) = \sum_{k=1}^{\infty} c^2(k) = ||f||^2_{\gamma}. \eqno(25)
 $$

 More exactly, we intend  to consistent estimate (consistently) the {\it truncated}
energy  functional of a view

 $$
 \underline{G} = \underline{G(f)} =
 \sum_{k=1}^{N_0(n)} c^2(k) = ||P_{N_0(n)}f||^2_{\gamma},\eqno(26)
 $$
where $ P_{N_0(n)}f $ denotes the projection of the function $ f(\cdot) $ onto
the subspace generated by the first $ N_0(n) $ polynomials $ \phi_k(\cdot), \
k=1,2,\ldots,. $ \par
 But it is evident that asymptotically as $ n \to \infty \
 \underline{G} \to G(f). $ \par
  If we define  as a first approximation the following estimation for the functional $ G $ the following expression

 $$
 G_0(n) = \sum_{k=1}^{M(n)} c^2(k,n),\eqno(27)
 $$
then we obtain

$$
G_0(n) \sim \sum_{k=1}^{M(n)}  (c(k) + \sigma K(\gamma) \eta(k)/\sqrt{n})^2,\eqno(28)
$$

$$
{\bf E} G_0(n) \sim \sum_{k=1}^{N_0(n)} c^2(k) + \sigma^2 K^2(\gamma) N_0(n)/n,\eqno(29)
$$
therefore an asymptotical non-biased  estimate of a $ G = G(f) $ based on the
$ G_0(n) $ has a view

$$
G(n): = \sum_{k=1}^{M(n)} c^2(k,n) - M(n) \ K^2(\gamma) \ \sigma^2(n)/n +
\tau^*(n). \eqno(30)
$$

{\bf Theorem 3.} As  $ n  \to \infty $ the r.v. $ G(n) $ has Gaussian distribution
with parameters

$$
G(n) \sim N(G, 4 \sigma^2 K^2(\gamma)G/n).\eqno(31)
$$

\vspace{3mm}
{\bf 2. Energy of general order.} Let us define

$$
W_w(f;n,N) = \sum_{k=1}^N  w(k) \left[c^2(k,n) \right] - n^{-1} K^2(\gamma) \sigma^2(n)+
\tau^*_w(n). \eqno(32)
$$

{\it We consider in this pilcrow only the case}
$ w(k) = w_{\theta}(k) \asymp k^{\theta},\theta = \const > 0,  \ k \to \infty; $ the
case when $ w(k) $ is regular varying as $ k \to \infty $ may be considered analogously.
\par
 We find by direct calculation {\it under assumption }
$  f \in W_{w^2}: $

 $$
 Z(W_w(f; n,N), W_w(f)) \asymp n^{-1} + \rho^2_w(N) + n^{-2} N^{2\theta+1}.\eqno(33)
 $$

 If we choose

 $$
 N_w = \argmin_N \left( \rho^2_w(f;N) + N^{2\theta+1}/n^2 \right),
 $$
we obtain the optimal in order {\it non-adaptive}  as $ n \to \infty $ estimation of
the energy $ W_w(f) $ of a view

$$
\hat{W}_w(f) = W_w(f; n, N_w).
$$

 For example, if $ \rho_w(f;N) \asymp N^{-2 \Delta}, \ \Delta = \const >0,  $  then

 $$
 N_w = \Ent[n^{1/(2\Delta + \theta + 1/2)}], \
 $$
$  \Ent[z] $ denotes the integer part of the number $ z, $ and

 $$
 Z(W_w(f; n, N_w),W_w(f)) \asymp
 \max \left(n^{-1}, n^{-4\Delta/(2\Delta + 2m+1/2 )}  \right).
 $$
 Notice that in the case when $ \Delta \ge 0.5 \theta + 1/4 $

$$
Z(W_w(f; n, N_w),W_w(f)) \asymp n^{-1}, \eqno(34)
$$
i.e. optimality in order as $ n \to \infty. $\par
\vspace{3mm}

 In order to construct an {\it adaptive } estimation of the variable $ W_w(f), $ we introduce the following statistics:

 $$
 \tau_w(N) = \tau_w(N,n) = \sum_{k=N+1}^{2N} w(k) \ c^2(k,n),
 $$

$$
\hat{N}_w = \argmin_{N \le n/3} \tau_w(N),
$$
and introduce the following {\it adaptive }estimation (measurement) of energy
$ W_w[f]: $
$$
\hat{W}_w = W_w(f; \hat{N}_w,n).
$$
\vspace{3mm}

{\bf Theorem 4.} Suppose $ f \in W_{2\theta}  $  and

 $$
\rho_w(f;N) \asymp N^{-2\Delta}L(N), \ \Delta = \const > 0,
 $$
where as before $ L = L(N) $ is slowly varying as $ N \to \infty $ function.
 Proposition:

 $$
 Z \left(\hat{W}_w, W_w \right) \le K_4(\Delta, \theta) \cdot
 \max \left(n^{-1}, n^{-4\Delta/(2\Delta +\theta + 1/2 )}  \right).\eqno(35)
 $$

 Notice that if $ \Delta \ge \theta/2 + 1/2, $ then

 $$
 Z \left(\hat{W}_w, W_w \right) \le K_5(\Delta, \theta) \cdot n^{-1},
 \eqno(36)
 $$
(optimality in order as $ n \to \infty $ of offered adaptive energy estimation).\par

\vspace{3mm}
\section{ Sketch of proofs.}\par
\vspace{3mm}

{\bf 1. Denoising.}  We have the expansion:

$$
c(k,n) = n^{-1} \sum_{i=1}^n \xi(i) \phi_k(x(i,n)) \gamma(x(i,n)) =
n^{-1} \sum_{i=1}^n \phi_k(x(i,n)) \gamma(x(i,n)) \ f(x(i,n)) +
$$

$$
n^{-1} \epsilon(i) \gamma(x(i,n)) \phi_k(x(i,n)) \sim c(k) +  n^{-1/2} \ \sigma \eta(k).
\eqno(37)
$$
 The r.v. $ \{ \eta(k) \} $ are centered, asymptotically independent and normal distributed  with variance as $ n \to \infty $

 $$
 {\bf Var} [\eta(k)] = n^{-1} \sum_{i=1}^n \phi^2_k(x(i,n)) \gamma^2(x(i,n)) \to
 \int_{-1}^1 \phi^2_k(x) \gamma^2(x) dx.
 $$
Further, as $ k \to \infty $

$$
{\bf Var}[ \eta(k)] \to \int_{-1}^1 \frac{\gamma(x)}{\sqrt{1-x^2}} dx = 2 \pi \ K(\gamma). \eqno(38)
$$
 On the other words, we get as a first approximation the following expression:

 $$
 c(k,n) \sim c(k) + C_1(\gamma) \ \sigma \ n^{-1/2} \zeta(k),
 $$
 where $ \{\zeta(k) \} $ are independent and standard normal distributed r.v. \par

We find by direct calculations:

$$
A(N,n) \asymp  \rho(N) + n^{-1} \ N; \eqno(39)
$$

$$
{\bf E}[\tau(N,n)] \asymp  \rho(N) + n^{-1} N = A(N,n);\eqno(40)
$$

$$
{\bf Var}[\tau(N,n)] \asymp n^{-2}N + n^{-1} \ \rho(N). \eqno(41)
$$

 We conclude by virtue of adaptive conditions that for some  positive constant
  $ C = C(\gamma) $

 $$
{\bf Var}[\tau(N,n)]  \le C \ n^{-1} \ \cdot A(N,n).
 $$
 Further considerations are alike to the \cite{Bobrov1}, \cite{Ostrovsky1}, chapter 5, section 13. \par
Notice that with probability one

$$
\min_{N \in [1,n]} \tau(N,n) \asymp  A^*(n), \ n \to \infty. \eqno(42)
$$
see \cite{Bobrov1}, \cite{Ostrovsky1}, chapter 5, section 13.\par
\vspace{2mm}
{\bf 2. Ordinary energy estimation.} We  have

$$
G(n) - G = \frac{2 \sigma K(\gamma)}{\sqrt{n}}\sum_{k=1}^{M(n)}
c(k)\eta(k) + \frac{\sigma^2 K^2(\gamma)}{n}\sum_{k=1}^{M(n)}(\eta^2(k)-1) \sim
$$

$$
\sim \frac{2 \sigma K(\gamma)}{\sqrt{n}}\sum_{k=1}^{N_0(n)}
c(k)\eta(k).
$$

 Therefore

 $$
 {\bf Var}[G(n) - G] \sim 4 \sigma^2 K^2(\gamma)n^{-1} \sum_{k=1}^{N_0(n)} c^2(k)
 \sim
 $$

 $$
 4 \sigma^2 K^2(\gamma)n^{-1} \sum_{k=1}^{\infty} c^2(k) = 4 \sigma^2 K^2(\gamma)G/n.
 \eqno(43)
 $$
 The asymptotical  normality follows from the classical CLT. \par
\vspace{3mm}
{\bf 3.} We prove here the implication $ (7a) \to (7b). $ \par
 Denote by $ F_n(x) $ the empirical function of distribution for the sequence
 $ x(i); \ F(x):= 0.5(x+1), \ x \in [-1,1]; $ obviously,

 $$
 \sup_{x \in [-1,1]} |F_n(x)-F(x)| \le 2/n.
 $$
 Further, put
 $$
 \sigma_{n,k} = \int_{-1}^1 [f(x)\phi_k(x)\gamma(x)] d(F_n(x) - F(x)),
 $$

 $$
 \Sigma(n) = ||f_{(N)} - f_{(N,n)}||^2_{\gamma};
 $$
then

$$
\Sigma(n) = \sum_{k=1}^{N_0} \sigma^2_{n,k}.
$$
We can and will assume without loss of generality that the function $ f(\cdot) $ has
a bounded derivative.\par
 We have for arbitrary values $ \epsilon \in (0,1/2): $

 $$
 |\sigma_{n,k}| \le  \int_{-1}^{-1+\epsilon} \left[f(x)\phi_k(x)\gamma(x)\right]
  d(F_n(x) - F(x)) +
 $$

$$
 \int_{-1+\epsilon}^{1-\epsilon} [f(x)\phi_k(x)\gamma(x)] d(F_n(x) - F(x)) +
$$

$$
 \int_{1-\epsilon}^{1} [f(x)\phi_k(x)\gamma(x)] d(F_n(x) - F(x)) = \sigma^{(1)}_{n,k} +
 \sigma^{(2)}_{n,k} + \sigma^{(3)}_{n,k};
$$

 $$
 |\sigma^{(1)}_{n,k}|\le C_1 \int_{-1}^{-1+\epsilon} \gamma(x) dx \le  C_2
 \epsilon^{1-\kappa}
 $$
and analogously may be estimated the value $ |\sigma^{(3)}_{n,k}|. $  \par
 Further,
$$
|\sigma^{(2)}_{n,k}| \le C_3 \left[\sigma^{(4)}_{n,k}+\sigma^{(5)}_{n,k} \right],
$$
where
$$
\sigma^{(4)}_{n,k} \stackrel{def}{=}
\int_{-1+\epsilon}^{1-\epsilon} |F_n(x)-F(x)||\phi^/_k(x)| \gamma(x) \ dx \le C_4k/n;
$$

$$
\sigma^{(5)}_{n,k} \stackrel{def}{=}
\int_{-1+\epsilon}^{1-\epsilon} |F_n(x)-F(x)||\phi_k(x)| |\gamma^/(x)| \ dx \le
C_5 \epsilon^{-\kappa}/n.
$$

Therefore,

$$
\sigma_{n,k} \le C_6 \left[ \frac{k}{n} + \frac{\epsilon^{-\kappa}}{n} + \epsilon^{1-\kappa}  \right].
$$
 We conclude choosing $ \epsilon = 0.25/n: $

$$
 \sigma^2_{n,k} \le C_7 \left[\frac{k^2}{n^2} + \frac{1}{n^{2(1-\kappa)}} \right];
$$

$$
\Sigma_{N_0} \le C_8 \left[ \frac{N_0^3}{n^2} + \frac{N_0}{n^{2(1-\kappa)}} \right]
$$
and following

$$
\lim_{n \to \infty} \frac{\Sigma_{N_0}}{A^*(n)}\le
\lim_{n \to \infty} \frac{\Sigma_{N_0}}{N_0/n}  = 0.
$$

\vspace{3mm}
{\bf 4. The proof of theorem 4} in completely analogous; see, for example,
\cite{Efromovich1} \cite{Klemela1}. Indeed, we have as $ n \to \infty $

$$
\hat{W_w}(n) - W_w \sim -\rho_W(N_w) +
\frac{2 \sigma K(\gamma)}{\sqrt{n}}\sum_{k=1}^{N_W(n)}
c(k)\eta(k)w(k) + \frac{\sigma^2 K^2(\gamma)}{n}\sum_{k=1}^{N_W(n)} w(k)(\eta^2(k)-1),
$$
therefore

$$
{\bf E}(\hat{W_w}(n) - W_w )^2 \sim \rho^2_W(N_w),
$$

$$
{\bf Var}(\hat{W_w}(n) - W_w ) \sim \frac{4\sigma^2 K^2(\gamma)}{n}
\sum_{k=1}^{N_w} c^2(k) w^2(k) +
\frac{3\sigma^4 K^4(\gamma)}{n^2} \sum_{k=1}^{N_w}  w^2(k);
$$
and we get under formulated above conditions, in particular, $ S_{w^2}(f) < \infty: $
$$
Z(\hat{W_w}(n),W_w ) \le C \left[\rho^2_W(N_w) + n^{-1} + S_{w^2}(f)/n^2 \right]
$$
etc.\par

\vspace{3mm}
\section{Concluding remarks, applications}
\vspace{3mm}

{\bf 0.} The optimality in order as $ n \to \infty $ of described above signal denoising
was proved by I.A.Ibragimov; see also the monograph \cite{Tchentsov1}, chapter 4. \\
{\bf 1.} We conclude making the so-called Fisher's transform for the ordinary energy estimation that

$$
\sqrt{G(n)} \sim N(\sqrt{G}, 4 \sigma^2 K^2(\gamma)/n).
$$
 This fact allows us to built a confidence interval for the value $ G. $ For example, with probability $ \approx 0.95 $

 $$
 |\sqrt{G} - \sqrt{G(n)}| \le \frac{6\sigma(n) K(\gamma)}{\sqrt{n}}.\eqno(44)
 $$
{\bf 2.} Let us return to the problem of diagnosis, for example, cardiodiagnostics.
The offered estimations with corresponding confidence regions allow us to detect the
distortions for obtained signal, for instance, in the cardiogram. \par
 Namely, if in the current time interval  the energy value does not belongs to the
 confidence region (or does not belongs to this region twice or three times,) we can
 conclude the presence of some distortion in the signal, in particular,
 infarction etc. \par

\end{document}